\documentclass[twoside,11pt]{amsart}
\usepackage{latexsym}
\usepackage{graphics, epsfig}
\newtheorem{theorem}{Theorem}[section]
\newtheorem{lemma}{Lemma}[section]
\newtheorem{proposition}{Proposition}[section]

\begin{document}

 \title{SIMPLY CONNECTED SYMPLECTIC $4$-MANIFOLDS \\
        WITH $b_{2}^+ =1$ AND $c_{1}^2 =2$}

 \author{Jongil Park}

 \address{Address 1: Department of Mathematics, Konkuk University,
          1 Hwayang-dong, Kwangjin-gu, Seoul 143-701, Korea}

 \address{Address 2: Department of Mathematics, Michigan
          State University, East Lansing, MI48824, USA}

 \email{jipark@konkuk.ac.kr\\ jipark@math.msu.edu}

 \thanks{This work was supported by grant No. R14-2002-007-01002-0
        from KOSEF}

 \date{Revised: March 23, 2004}



\begin{abstract}
 In this article we construct a new family of simply connected
 symplectic $4$-manifolds with $b_{2}^+ =1$ and $c_{1}^2 =2$
 which are not diffeomorphic to rational surfaces by using
 rational blow-down technique.
 As a corollary, we conclude that a rational surface
 ${\mathbf CP}^2 \sharp 7{\overline{{\mathbf CP}}^2}$
 admits an exotic smooth structure.
\end{abstract}

\maketitle

\section{Introduction}

\markboth{JONGIL PARK}{SIMPLY CONNECTED SYMPLECTIC $4$-MANIFOLDS
                       WITH $b_{2}^+ =1$ AND $c_{1}^2 =2$}

 One of the fundamental problems in the topology of $4$-manifolds is to
 determine whether a given topological $4$-manifold admits a
 smooth structure and, if it does, whether such a smooth structure
 is unique or not. Though the complete answer is far from reach,
 gauge theory makes us to answer partially these questions
 (\cite{fs1}, \cite{fs2}, \cite{g}, \cite{sz}, \cite{t}).
 But most known results are in the case of simply connected $4$-manifolds
 with either \mbox{$b_{2}^+ > 1$} odd or $b_{2}^+ =1$ and $c_{1}^2 \leq 0$.
 In the case when $b_{2}^+ =1$ and $c_{1}^2 > 0$, a theorem of D. Kotschick
 is the only known result that the Barlow surface is not
 diffeomorphic to
 ${\mathbf CP}^2 \sharp 8{\overline{{\mathbf CP}}^2}$ (\cite{k1}).
 Since then, there was little progress on the problems.

 In this paper we investigate exotic smooth structures on a
 rational surface
 ${\mathbf CP}^2 \sharp 7{\overline{{\mathbf CP}}^2}.$
 According to a convention, we say
 that a smooth $4$-manifold admits an {\em exotic smooth
 structure} if it has more than one distinct smooth structure.
 One way to get an exotic smooth structure on a given smooth
 $4$-manifold $X$ is to construct a new smooth $4$-manifold $X'$
 which is homeomorphic, but not diffeomorphic, to $X$.
 Hence the problem of finding exotic smooth structures on a rational
 surface ${\mathbf CP}^2 \sharp 7{\overline{{\mathbf CP}}^2}$ is
 equivalent to find a new family of simply connected
 smooth $4$-manifolds with $b_{2}^+ =1$ and $c_{1}^2 =2$.
 Note that the only known simply connected closed smooth $4$-manifolds
 with $b_{2}^+ =1$ and $c_{1}^2 \geq 2$ are {\mbox rational} surfaces such
 as ${\bf CP}^2, S^2\times S^2$ and ${\mathbf CP}^2 \sharp
 n{\overline{{\mathbf CP}}^2} (n \leq 7)$.
 Despite the fact that it is no constraint on the existence of
 simply connected smooth \mbox{$4$-manifolds} with $b_{2}^+ =1$
 and $c_{1}^2 \geq 2$ which are not rational surfaces, no such $4$-manifolds
 have been known. Thus it has long been an interesting question to find such
 $4$-manifolds (refer to Problem 4.45 in the Kirby list appeared
 in \cite{ka}).
 In \mbox{Section 3} we construct a new family of simply connected
 symplectic $4$-manifolds  with $b_{2}^+ =1$ and $c_{1}^2 =2$
 which are homeomorphic, but not diffeomorphic, to rational surfaces
 by using rational blow-down technique introduced by R. Fintushel and
 R. Stern in~\cite{fs1}. As one of our main results, we get the following

\begin{theorem}
\label{thm-main}
 There exists a simply connected symplectic $4$-manifold with
 $b_{2}^+ =1$ and $c_{1}^2 =2$ which is homeomorphic, but not
 diffeomorphic, to a rational surface
 ${\mathbf CP}^2 \sharp 7{\overline{{\mathbf CP}}^2}$.
\end{theorem}

 Hence we conclude that a rational surface
 ${\mathbf CP}^2 \sharp 7{\overline{{\mathbf CP}}^2}$
 admits an exotic smooth structure.
 Furthermore, by blowing up of a symplectic $4$-manifold
 constructed in Theorem~\ref{thm-main} above,
 we also conclude that a rational surface
 ${\mathbf CP}^2 \sharp 8{\overline{{\mathbf CP}}^2}$
 admits at least three distinct smooth structures -
 an Einstein metric with positive scalar curvature,
 an Einstein metric with negative scalar curvature and
 no Einstein metric.
 \\

 {\em Acknowledgements}. The author would like to thank Ronald Fintushel
 for his constant encouragement and guidance. The author is deeply
 indebted to him for sharing his ideas with the author while
 working on this problem. The author would also like to thank Andr\'{a}s
 Stipsicz for pointing out some critical errors in the earlier
 version of this article and an explanation about various types of singular
 fibers appeared in a rational surface $E(1)$. The author also
 thanks Denis Auroux for figuring out $\widetilde{E_6}$-fiber in $E(1)$
 and Tian-Jun Li for an explanation about the symplectic $2$-forms on a
 non-minimal rational surface. Finally he also thanks the department of
 Mathematics at Michigan State University for hospitality and support
 during his visit. \\

\section{Preliminaries}
\label{sec-2}
 In this section we briefly review the Seiberg-Witten theory and a
 rational blow-down surgery which will be the main technical tools
 to get our results.

 First we briefly introduce the Seiberg-Witten theory for smooth
 $4$-manifolds. In particular, we pay attention to the Seiberg-Witten
 invariant of $4$-manifolds with $b_{2}^+ =1$ (\cite{m} for details).
 Let $X$ be a closed, oriented smooth $4$-manifold with $b_{2}^+ >0$
 and a fixed metric $g$, and let $L$ be a characteristic line bundle on $X$,
 i.e. $c_{1}(L)$ is an integral lift of  $w_{2}(X)$. This determines a
 $Spin^{c}$-structure on $X$ which induces a complex spinor
 bundle $W \cong W^{+} \oplus W^{-}$, where $W^{\pm}$ is
 the associated  $U(2)$-bundles on $X$ such that
 $\mathrm{det}(W^{\pm}) \cong L$.
 Note that the Levi-Civita connection on $TX$  together with a unitary
 connection $A$  on $L$ induces a connection  $\nabla_{A} :
 \Gamma(W^{+}) \rightarrow \Gamma(T^{\ast}X\otimes  W^{+})$.
 This connection,  followed  by  Clifford multiplication, induces a
 $Spin^{c}$-Dirac operator
 $D_{A} :  \Gamma(W^{+}) \rightarrow \Gamma(W^{-})$.
 Then, for each self-dual $2$-form
 $h \in \Omega^2_{+_{g}}(X\!:\!{\mathbf R})$,
 the following pair of equations  for a unitary connection $A$ on $L$
 and a section $\Psi$ of $\Gamma (W^{+})$ are called the
 {\em perturbed Seiberg-Witten equations}:
\begin{equation}
  (SW_{g,h})\left\{ \begin{array}{ll}
              D_{A}\Psi  & =\ \  0     \label{sw-eq}\\
              F_{A}^{+_{g}} & = \ \ i(\Psi \otimes \Psi^{\ast})_{0} + ih\, .
                  \end{array} \right.
\end{equation}
 Here $F_{A}^{+_{g}}$ is the self-dual part of the curvature of $A$ with
 respect to a metric $g$ on $X$ and $(\Psi  \otimes \Psi^{\ast})_{0}$ is the
 trace-free part of $(\Psi \otimes  \Psi^{\ast})$ which is interpreted as an
 endomorphism of $W^{+}$. The gauge group $\mathcal{G}:=Aut(L)\cong
 Map(X,S^{1})$ acts on the  space
 $\mathcal{A}_{X}(L) \times \Gamma (W^{+})$  by
 \[g\cdot (A,\Psi) =(g\circ A\circ g^{-1},g\cdot \Psi) \]
 Since the set of solutions is invariant under the action, it induces
 an orbit space, called the {\em Seiberg-Witten  moduli space}, denoted by
 $M_{X,g,h}(L)$, whose formal dimension is
 \[ \mathrm{dim}M_{X,g,h}(L)=\frac{1}{4}(c_{1}(L)^{2}-3\sigma(X)-2e(X))\]
 where $\sigma(X)$ is the signature of $X$ and $e(X)$ is the Euler
 characteristic of  $X$. Note that if $b_{2}^{+}(X) > 0$ and
 $M_{X,g,h}(L)\neq \phi$, then for a generic self-dual $2$-form $h$ on $X$
 the moduli space $M_{X,g,h}(L)$ contains no reducible solutions,
 so that it is a compact, oriented, smooth manifold of the given dimension.\\

\noindent
 {\bf Definition} The {\em Seiberg-Witten invariant}
 ({\em for brevity, SW-invariant}) for a smooth $4$-manifold $X$ with
 $b_{2}^+ >0$ is a function
 $SW_{X} : Spin^{c}(X) \rightarrow  {\mathbf Z}$  defined by
 \begin{equation}
 SW_{X}(L) := \left\{ \begin{array}{ll}
  <\beta^{d_{L}},[M_{X,g,h}]>  &  \mathrm{if\ \ dim}M_{X,g,h}(L):=2d_{L}
    \geq 0\ \ \mathrm{and\ \ even} \\
  \ \ \ \ \ \ \ 0 &  \mathrm{otherwise} \, .
                       \end{array}
               \right.
\end{equation}
 Here $\beta$ is a generator of $H^{2}(\mathcal{B}_{X}^*(L);{\bf
 Z})$ which is the first Chern class of the $S^1$-bundle
 \[ \widetilde{\mathcal{B}}_{X}^*(L) =
    \mathcal{A}_{X}(L) \times (\Gamma (W^{+})\!-\!\{0\})/Aut^{0}(L)
    \longrightarrow \mathcal{B}_{X}^*(L) \]
 where $Aut^{0}(L)$ consists of gauge transformations which are
 the identity on the fiber of $L$ over a fixed base point in $X$.
 Note that if $b_{2}^{+}(X) >1$, the Seiberg-Witten invariant, denoted by
 $SW_{X}= \sum SW_{X}(L)\cdot e^{c_{1}(L)}$, is a diffeomorphism invariant,
 i.e. $SW_{X}$ does not depend on the choice of a metric on $X$ and a
 generic perturbation of Seiberg-Witten equations.
 Furthermore, only finitely many $Spin^{c}$-structures on $X$ have
 a non-zero Seiberg-Witten  invariant. We say that the characteristic
 line bundle $L$, equivalently a  cohomology class
 $c_{1}(L) \in H^{2}(X;{\mathbf Z})$,
 is a {\em SW-basic class} of $X$ if $SW_{X}(L)\neq 0$.

 When $b_{2}^{+}(X)\! =\! 1$, the SW-invariant $SW_{X}(L)$ defined
 in $(2)$ above depends not only on a metric $g$ but also on a self-dual
 $2$-form $h$. Because of this fact, there are several types of
 Seiberg-Witten invariants for a smooth $4$-manifold with $b_{2}^+\! =\! 1$
 depending on how to perturb the Seiberg-Witten equations.
 We introduce two types of SW-invariants and investigate how they
 are related. First we allow all metrics and self-dual $2$-forms to perturb
 the Seiberg-Witten equations. Then the SW-invariant $SW_{X}(L)$ defined
 in $(2)$ above has generically two values which are determined by the
 sign of $(2\pi c_{1}(L) +[h])\cdot [\omega_{g}]$, where $\omega_{g}$
 is a  unique $g$-self-dual harmonic $2$-form of norm one lying in the
 (preassigned) positive component of $H^{2}_{+_{g}}(X;{\mathbf R})$.
 We denote the SW-invariant for a metric $g$ and a generic self-dual $2$-form
 $h$ satisfying $(2\pi c_{1}(L) +[h])\cdot [\omega_{g}]>0$ by $SW_{X}^+(L)$
 and denote the other one by $SW_{X}^-(L)$. Then the wall crossing
 formula tells us the relation between $SW_{X}^+(L)$ and $SW_{X}^-(L)$.

\begin{theorem}[Wall crossing formula, \cite{m}]
\label{wall}
 Suppose that $X$ is a closed, oriented smooth $4$-manifold with
 $b_{1}=0$ and $b_{2}^+ =1$. Then for each characteristic line bundle
 $L$ on $X$ such that the formal dimension of the moduli space
 $M_{X,g,h}(L)$ is non-negative and even, say $2d_{L}$, we have
 \[ SW_{X}^+(L)  -  SW_{X}^-(L)  = -(-1)^{d_{L}} \, .\]
\end{theorem}

\vspace{.2cm}

  By the way, C. Taubes' result on the SW-invariant of a symplectic
 $4$-manifold with $b_{2}^+ >1$ can be easily extended to the
 $b_{2}^+ =1$ case.

\begin{theorem}[\cite{t}, \cite{ll1}]
\label{T-2}
 Suppose $X$ is a closed symplectic $4$-manifold with $b_{2}^+ =1$
 and a canonical class $K_{X}$. Then $SW_{X}^{-}(-K_{X}) = \pm 1$.
\end{theorem}

\vspace{.2cm}

  Second one may perturb the Seiberg-Witten equations by adding only a
 small generic self-dual $2$-form $h \in \Omega^2_{+_{g}}(X; {\mathbf R})$,
 so that one can define the SW-invariants as in $(2)$ above. In this case
 we denote the SW-invariant for a metric $g$ satisfying
 $(2\pi c_{1}(L))\cdot [\omega_{g}]>0$ by $SW_{X}^{\circ,+}(L)$
 and we denote the other one by  $SW_{X}^{\circ,-}(L)$.
 Note that $SW_{X}^{\circ,\pm}(L)=SW_{X}^{\pm}(L)$.
 But it sometimes happens that the sign of
 $(2\pi c_{1}(L))\cdot [\omega_{g}]$ is the same for all generic metrics,
 so that there exists only one SW-invariant obtained by a small generic
 perturbation of the Seiberg-Witten equations.
 In such a case we define the SW-invariant of $L$ on $X$ by
\begin{equation*}
 SW_{X}^{\circ}(L) := \left\{ \begin{array}{ll}
 SW_{X}^{\circ, +}(L)  &  \mathrm{if}\  2\pi c_{1}(L) \cdot [\omega_{g}] > 0 \\
 SW_{X}^{\circ, -}(L)  &   \mathrm{if}\  2\pi c_{1}(L) \cdot [\omega_{g}] < 0 \, .
                       \end{array}
               \right.
\end{equation*}
 If $SW_{X}^{\circ}(L) \neq 0$, we call the corresponding
 $c_{1}(L)$ (or $L$) a {\em SW-basic class} of $X$.
 Then the Seiberg-Witten invariant of $X$, denoted by
 $SW_{X}^{\circ}= \sum SW_{X}^{\circ}(L)\cdot e^{c_1(L)}$,
 will also be a diffeomorphism invariant.
 Furthermore we can extend many results obtained for smooth
 $4$-manifolds with $b_{2}^+ >1$ to this case.
 For example, we have

\begin{theorem}
\label{thm-basic}
 Let $X$ be a  simply connected closed smooth $4$-manifold with
 $b_{2}^+ =1$ and $b_{2}^- \leq 9$. Then  \\
 $($i$)$ There are only finitely many characteristic line bundles $L$ on $X$
         such that \mbox{$SW_{X}^{\circ}(L) \neq 0$}. \\
 $($ii$)$ If $X$ admits a metric of positive scalar curvature, then the
      SW-invariant of $X$ vanishes, that is, $SW_{X}^{\circ}(L) = 0$ for any
      characteristic line bundle $L$ on $X$.
\end{theorem}

 {\em Proof} : Proofs of (i) and (ii) are exactly the same as the
 proofs of case $b_{2}^+ >1$ as long as the SW-invariant $SW_{X}^{\circ}$
 is well defined,  i.e. it is independent of metrics on $X$.
 Let $L$ be a characteristic line  bundle on $X$ such that the formal
 dimension, $\frac{1}{4}(c_{1}(L)^2 - 3\sigma(X) -2e(X))$,
 of the moduli space is non-negative and even.
 The condition  $b_{2}^+ =1$ and $b_{2}^- \leq 9$
 imply  that  $c_{1}(L)^2 \geq 3\sigma(X) +2e(X) \geq 0$.
 Furthermore, since $X$ is simply  connected and $c_{1}(L)$ is
 characteristic,  $c_{1}(L) \neq 0$ (Otherwise, $X$ has $b_{2}^{-}=9$
 and it is spin which contradicts the Rohlin's signature theorem.)
 Thus, for any metric $g$ on $X$, the light cone lemma implies
 $c_{1}(L)\cdot [\omega_{g}] \neq 0$, so that the sign of
 $(2\pi c_{1}(L))\cdot [\omega_{g}]$ is the same for all
 generic metrics. Hence the SW-invariant $SW_{X}^{\circ}(L)$
 is well defined. $\ \ \ \Box$ \\

 Next we briefly review a {\em rational blow-down} technique
 introduced by R. Fintushel and R. Stern and state related facts
 (\cite{fs1} for details).

 Let $C_{p}$ be a smooth $4$-manifold obtained by plumbing the
 $(p-1)$ disk bundles over the $2$-sphere instructed by
 the following diagram

 \begin{picture}(400,60)(-100,-25)
   \put(-12,3){\makebox(200,20)[bl]{$-(p+2)$ \hspace{4pt}
                                    $-2$ \hspace{92pt} $-2$}}
   \put(4,-25){\makebox(200,20)[tl]{$u_{p-1}$ \hspace{10pt}
                                    $u_{p-2}$ \hspace{86pt} $u_{1}$}}
   \multiput(10,0)(40,0){2}{\line(1,0){40}}
   \multiput(10,0)(40,0){2}{\circle*{3}}
   \multiput(100,0)(5,0){4}{\makebox(0,0){$\cdots$}}
   \put(125,0){\line(1,0){40}}
   \put(165,0){\circle*{3}}
 \end{picture}

 \noindent
 where each vertex $u_{i}$ represents a disk bundle over
 the $2$-sphere with Euler class labelled above and an interval
 between vertices indicates plumbing the disk bundles corresponding
 to the vertices. Label the homology classes represented by the
 $2$-spheres in $C_{p}$ by $u_{1}, \ldots, u_{p-1}$ so that the
 self-intersections are $u_{p-1}^2 =-(p+2)$ and $u_{i}^2= -2$ for
 $1 \leq i \leq p-2$. Furthermore, orient the $2$-spheres so that
 $u_{i}\cdot u_{i+1} = +1$. Then a configuration $C_{p}$ has
 the following topological properties:\\

\begin{enumerate}
 \item It is a negative definite simply connected smooth $4$-manifold
       whose boundary is the lens space $L(p^2,1-p)$, and
       the lens space $L(p^2,1-p) = \partial{C_{p}}$ bounds a rational
       ball $B_{p}$ with $\pi_{1}(B_{p}) \cong {\mathbf Z}_{p}$.
 \item $H_{2}(C_{p};{\mathbf Z}) \cong \bigoplus_{i=1}^{p-1}{\mathbf Z}$
       has generators $\{u_{i} :1 \leq i \leq p-1\}$, where each $u_{i}$
       can be represented by the zero-section $S_{i}^{2}$ of the disk bundle
        $u_{i}$ over $S^{2}$ (We use $u_{i}$ for both a generator and
        the corresponding disk bundle).
 \item  Let $P$ be a plumbing matrix for $C_{p}$ with respect to
       the basis $\{u_{i} :1 \leq i \leq p-1\}$.
       Then the intersection form on $H^{2}(C_{p};{\mathbf Q})$
       with respect to the dual basis $\{\gamma_{i}: 1\leq i\leq p-1\}$
       (i.e. $<\gamma_{i}\,\,,\,u_{j}>=\delta_{ij}$)  is given by
       \[ Q := (\gamma_{i}\cdot \gamma_{j}) = P^{-1} . \]
       For example, when $p=7$, we have the following intersection form
       on $H^{2}(C_{7};{\mathbf Q})$ with respect to
       $\{\gamma_{i}: 1\leq i\leq 6\}$:
       \[ Q_7 = \frac{-1}{49} \left( \begin{array}{cccccc}
                        41 & 33 & 25 & 17 &  9 & 1  \\
                        33 & 66 & 50 & 34 & 18 & 2  \\
                        25 & 50 & 75 & 51 & 27 & 3  \\
                        17 & 34 & 51 & 68 & 36 & 4  \\
                        9  & 18 & 27 & 36 & 45 & 5  \\
                        1  &  2 &  3 &  4 &  5 & 6
                                            \end{array}
                                            \right)  \]
\\
\end{enumerate}

\noindent
 {\bf Definition} Suppose $X$ is a smooth $4$-manifold which contains
 a configuration $C_{p}$. Then we construct a new smooth $4$-manifold
 $X_{p}$, called the {\em rational blow-down} of $X$, by replacing
 $C_{p}$ with a rational ball $B_{p}$.
 We call this a {\em rational blow-down} technique.
 Note that this process is well-defined, that is, a new smooth
 $4$-manifold $X_{p}$ is uniquely constructed (up to diffeomorphism)
 from $X$ because each diffeomorphism of $\partial B_{p}\! =\!
 L(p^2, 1\!-\!p)$ extends over the rational ball $B_{p}$.
 Furthermore,
 M. Symington proved that a rational blow-down manifold
 $X_{p}$ admits a symplectic structure in some cases.

\begin{theorem}[\cite{sy}]
\label{thm-sy}
 Suppose $X$ is a symplectic $4$-manifold containing a configuration
 $C_{p}$ with a symplectic $2$-form $\omega$. If all $2$-spheres
 $u_{i}$ in $C_{p}$ are symplectically embedded and intersect positively,
 then the rational blow-down manifold $X_{p}= X_{0}\cup_{L(p^2,1-p)}B_{p}$
 admits a symplectic $2$-form $\omega_p$ such that $(X_{0},
 \omega_p|_{X_{0}})$ is symplectomorphic to $(X_{0}, \omega|_{X_{0}})$.
\end{theorem}

 {\em Remark 1.} Suppose $X=X_{0}\cup_{L(p^2,1-p)}C_{p}$ is a symplectic
 $4$-manifold with a canonical class $K$ and a compatible
 symplectic $2$-form $\omega$. In the case when
 $X_{p}=X_{0}\cup_{L(p^2,1-p)}B_{p}$ admits a symplectic structure
 as in the Theorem~\ref{thm-sy} above, let $\omega_{p}$ be the induced
 symplectic $2$-form on $X_{p}$ such that
 $\psi_p:(X_{0}, \omega_p|_{X_{0}}) \rightarrow (X_{0}, \omega|_{X_{0}})$
 is a symplectomorphism. We also let $K_p$ be the canonical class on $X_p$
 which is induced from the symplectic $2$-form $\omega_p$ on $X_p$.
 Then, since $H^1(L(p^2,1-p);{\bf Q})=H^2(L(p^2,1-p);{\bf Q})=0$,
 if we decompose \mbox{$K$ and $\omega$ as}
 \[K = K|_{X_{0}} + K|_{C_{p}} \ \ \mathrm{and}\ \
   [\omega] = [\omega|_{X_{0}}] + [\omega|_{C_{p}}] \]
 with $K|_{X_{0}}, [\omega|_{X_{0}}] \in H^2(X_{0};{\bf Q})$ and
      $K|_{C_{p}}, [\omega|_{C_{p}}] \in H^2(C_{p};{\bf Q})$,
 we can also decompose \mbox{$K_p$ and $\omega_{p}$ as}
 \[ K_p = K_p|_{X_{0}} + K_p|_{B_{p}} \ \ \mathrm{and}\ \
   [\omega_{p}] = [\omega_{p}|_{X_{0}}] + [\omega_{p}|_{B_{p}}]. \]

\vspace{.2cm}

\begin{lemma}
\label{lem-can}
 Under the same hypothesis on $(X, K, \omega)$ and
 $(X_{p}, K_p, \omega_{p})$ as above, we have
 \[K_p \cdot [\omega_{p}]
   = K \cdot [\omega] - K|_{C_{p}} \cdot [\omega|_{C_{p}}]. \]
\end{lemma}

 {\em Proof} : Since $K_{p}|_{B_{p}}$ and $[\omega_{p}|_{B_p}]$ are
 zero elements in $H^2(B_{p};{\bf Q})$, we have
\begin{eqnarray*}
  K_p \cdot [\omega_{p}]
  & = & K_p|_{X_{0}} \cdot [\omega_p|_{X_{0}}]
    =  \psi_p^*(K|_{X_{0}}) \cdot \psi_p^*([\omega|_{X_{0}}]) \\
  & = & K|_{X_{0}} \cdot [\omega|_{X_{0}}]
    =   K \cdot [\omega] - K|_{C_{p}} \cdot [\omega|_{C_{p}}].
  \ \ \ \ \ \ \ \ \Box \\
\end{eqnarray*}

\section{A Main Construction}
\label{sec-3}
 In this section we construct a new family of simply connected
 symplectic $4$-manifolds with $b_{2}^+ =1$ and $c_{1}^2 =2$ which are
 homeomorphic, but not diffeomorphic, to rational surfaces by using rational
 blow-down technique introduced by R. Fintushel and R. Stern in~\cite{fs1}.

  Let us start with analyzing a simply connected rational surface
 $E(1)={\mathbf CP}^2 \sharp 9 {\overline{{\mathbf CP}}^2}$.
 There are several ways to describe $E(1)$. One way to construct $E(1)$
 is to take two generic cubic curves in ${\mathbf CP}^2$ which intersect
 each other at $9$ points and then blow up $9$ times at these points in
 ${\mathbf CP}^2$. This viewpoint makes us to see
 $E(1)={\mathbf CP}^2 \sharp 9 {\overline{{\mathbf CP}}^2}$
 as a Lefschetz fibration over ${\mathbf CP}^1$ whose generic fiber
 is an elliptic curve, say $f$, and which also has $6$ singular cusp fibers
 (or, equivalently, $12$ singular fishtail fibers).
 Since $4$ singular cusp fibers in $E(1)$ can be deformed to an
 $\widetilde{E_6}$-singular fiber, $E(1)$ can also be described as an
 elliptic fibration over ${\mathbf CP}^1$ with $3$ singular fibers,
 one $\widetilde{E_6}$-fiber and two cusp fibers. Note that a
 neighborhood of the $\widetilde{E_6}$-fiber in $E(1)$ is a smooth
 $4$-manifold obtained by plumbing disk bundles over the holomorphically
 embedded $2$-spheres $S_i (1\leq i \leq 7)$ of square $-2$ instructed
 by the Dynkin diagram of $\widetilde{E_6}$ (\cite{hkk} for details).

 \begin{figure}[hbtp]
 \begin{picture}(400,90)(-110,-10)
   \put(2,3){\makebox(200,20)[bl]{$-2$ \hspace{21pt} $-2$ \hspace{27pt}
               $-2$ \hspace{10pt} $-2$ \hspace{18pt} $-2$}}
   \put(4,-25){\makebox(200,20)[tl]{$S_{1}$ \hspace{22pt} $S_{2}$ \hspace{22pt}
              $S_{3}$ \hspace{22pt} $S_{4}$ \hspace{22pt} $S_{5}$ \hspace{22pt}}}
   \put(75,58){$S_{7}$  \hspace{0pt} $-2$}
   \put(75,28){$S_{6}$  \hspace{0pt} $-2$}
   \multiput(90,60)(0,-30){2}{\line(0,-1){30}}
   \multiput(90,60)(0,-30){2}{\circle*{3}}
   \multiput(10,0)(40,0){4}{\line(1,0){40}}
   \multiput(10,0)(40,0){5}{\circle*{3}}
 \end{picture}
 \caption{$\widetilde{E_6}$-singular fiber}
 \label{E6}
 \end{figure}
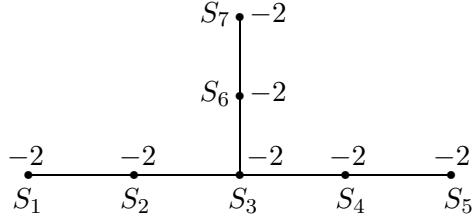

 \begin{lemma}[\cite{a}, \cite{f}]
 \label{lem-E6}
 The second (co)homology classes $[S_i]\, \, (1\leq i \leq 7)$ of the
 $2$-spheres $S_i$ embedded in $\widetilde{E_6}$ can be represented by
 $[S_1] =  e_4 - e_7,\, [S_2] =  e_1 - e_4,\, [S_3] = h -e_1 - e_2 -
 e_3,$  $[S_4] =  e_2 - e_5,\, [S_5] =  e_5 - e_9,\, [S_6] =  e_3 - e_6$
 and $[S_7] =  e_6 - e_8$, where $h$ denotes a generator of
 $H_{2}({\mathbf CP}^2;{\mathbf Z})$ and each $e_{i}$ denotes the
 (co)homology class represented by the $i^{th}$ exceptional curve in
 $\overline{{\mathbf CP}}^2 \subset
 E(1)={\mathbf CP}^2 \sharp 9 {\overline{{\mathbf CP}}^2}$.
 \end{lemma}

 {\em Proof} : Note that $E(1)$ can be constructed as follows:
 First choose a generic cubic curve $C$ (represented homologically
 by $3h$) which intersects with a line (represented by $h$) at $3$
 points in ${\mathbf CP}^2$. And then blow up at these $3$ points,
 so that we get an embedded $2$-sphere $S_3$, represented by
 $h-e_1-e_2-e_3$, of multiplicity $3$ and of square $-2$ in
 ${\mathbf CP}^2 \sharp 3{\overline{{\mathbf CP}}^2}$. Again blow
 up $3$ times at the intersection points between the curve
 $C-e_1-e_2-e_3$ and $3$ exceptional curves $e_1, e_2$ and $e_3$
 respectively, so that we get embedded $2$-spheres $S_2=e_1-e_4,
 S_4=e_2-e_5$ and $S_6=e_3-e_6$ of multiplicities $2$ and of squares
 $-2$ in ${\mathbf CP}^2 \sharp 6{\overline{{\mathbf CP}}^2}$.
 Finally blow up $3$ times at the intersection points
 between the curve $C-e_1-e_2- \cdots -e_6$ and $3$ new exceptional
 curves $e_4, e_5$ and $e_6$ respectively, so that
 we get again embedded $2$-spheres $S_1=e_4-e_7, S_5=e_5-e_9$ and
 $S_7=e_6-e_8$ of multiplicities $1$ and of squares
 $-2$ in ${\mathbf CP}^2 \sharp 9{\overline{{\mathbf CP}}^2}$.
 Then the embedded $2$-spheres $\{S_1\, \ldots \, S_7\}$ consists
 of $\widetilde{E_6}$-singular fiber in $E(1)$. $\ \ \ \Box$ \\

 Note that the standard canonical class
 $K_{E(1)}\in H^{2}(E(1);{\mathbf Z})$ of $E(1)$ is represented by
 $K_{E(1)} = -3h+(e_{1}+ \cdots + e_{9}) = -[f]$.
 Furthermore, there is a relation between the canonical class and
 a compatible symplectic $2$-form on a non-minimal rational surface
 which will play an important role in the proof of our main results.
 Explicitly, we have

\begin{lemma}
\label{lem-E(1,k)}
 For each integer $k \geq 1$, there exists a symplectic $2$-form
 $\omega$ on $E(1)\sharp k{\overline{\mathbf CP}}^2$
 which is compatible with the standard canonical class
 $K_{E(1)\sharp k{\overline{\mathbf CP}}^2}=-3h+(e_{1}+ \cdots + e_{9+k})$
 such that its cohomology class $[\omega]$ can be represented by
 $ah - (b_{1}e_{1} + \cdots + b_{9+k}e_{9+k})$ for some rational
 numbers $a, b_{1}, \ldots, b_{9+k}$ satisfying
 $a \geq b_1 \geq b_2 \geq \ldots \geq b_{9+k} \geq 0$ and
 $3a > b_{1} + \cdots + b_{9+k}$.
\end{lemma}

 {\em Proof} : Since $E(1)\sharp k{\overline{\mathbf CP}}^2$ is
 a rational surface, there exists a symplectic $2$-form
 $\omega$ on $E(1)\sharp k{\overline{\mathbf CP}}^2$
 which is compatible with the standard canonical class
 $K_{E(1)\sharp k{\overline{\mathbf CP}}^2}=-3h+(e_{1}+\cdots+e_{9+k})$
 satisfying $K_{E(1)\sharp k{\overline{\mathbf CP}}^2} \cdot [\omega] < 0$
 (refer to Corollary 1.4 in~\cite{ms}).
 Furthermore, Lemma 4.7 in~\cite{ll2} guarantees that the cohomology
 class $[\omega]$ of the symplectic $2$-form $\omega$ can be
 represented by $ah - (b_{1}e_{1} + \cdots + b_{9+k}e_{9+k})$
 for some rational numbers $a, b_{1}, \ldots, b_{9+k}$ satisfying
 $a \geq b_1 \geq b_2 \geq \ldots \geq b_{9+k} \geq 0$.
 The inequality $3a > b_{1} + \cdots + b_{9+k}$ follows from the fact
 that $K_{E(1)\sharp k{\overline{\mathbf CP}}^2} \cdot [\omega] < 0$.
 $\ \ \ \Box$ \\

 In fact, T. Li and A. Liu obtained many results regarding
 symplectic structures and canonical classes on
 symplectic $4$-manifolds with $b_2^+ =1$. For example, they
 proved

\begin{lemma}[\cite{ll1}]
\label{lem-ll1}
 There is a unique symplectic structure on
 ${\mathbf CP}^2 \sharp k{\overline{\mathbf CP}}^2$
 for $2\leq k \leq 9$ up to diffeomorphisms and deformation.
 For $k\geq 10$, the symplectic structure is still unique
 for the standard canonical class. In particular,
 ${\mathbf CP}^2 \sharp k{\overline{\mathbf CP}}^2$ $(2\leq k \leq
 9)$ does not admit a symplectic $2$-form $\omega$ for which
 $c_1(K)\cdot [\omega] >0$.
\end{lemma}

\begin{proposition}
\label{pro-C7}
 There exists a configuration $C_7$ in a rational surface
 $E(1)\sharp 4{\overline{\mathbf CP}}^2$ such that all $2$-spheres $u_{i}$
 lying in $C_7$ are symplectically embedded.
\end{proposition}

\begin{figure}[h]
 \begin{picture}(400,30)(-110,-5)
   \put(-15,3){\makebox(200,20)[bl]{\hspace{15pt} $-9$ \hspace{18pt} $-2$
                 \hspace{95pt} $-2$}}
   \put(-10,-25){\makebox(200,20)[tl]{\hspace{12pt} $u_{6}$ \hspace{22pt}
                 $u_{5}$ \hspace{98pt} $u_{1}$}}
   \multiput(10,0)(40,0){2}{\line(1,0){40}}
   \multiput(10,0)(40,0){2}{\circle*{3}}
   \multiput(100,0)(5,0){4}{\makebox(0,0){$\cdots$}}
   \put(125,0){\line(1,0){40}}
   \put(165,0){\circle*{3}}
 \end{picture}
 \caption{$C_7 \subset E(1)\sharp 4{\overline{\bf CP}}^2$}
 \label{C7}
 \end{figure}
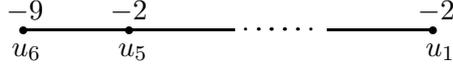

 {\em Proof} : Note that $E(1)$ can be viewed as an elliptic fibration
 with an $\widetilde{E_6}$-singular fiber and $2$ cusp
 singular fibers (equivalently, $4$ singular fishtail fibers).
 Since the homology class $[f]$ of the elliptic fiber $f$ in $E(1)$
 can be represented by an immersed $2$-sphere with one positive double
 point (equivalently, a fishtail fiber) and since $E(1)$ contains at least
 \mbox{$4$ such} immersed $2$-spheres, we blow up $4$ times at these double
 points so that there exist embedded $2$-spheres,
 $f-2e_{10}, \, \ldots, f-2e_{13}$, in $E(1)\sharp 4{\overline{\bf CP}}^2$
 which intersect a section $e_9$ of $E(1)$ positively at points,
 say $p_{1}, \, \ldots, p_{4}$, respectively.
 And then, resolving symplectically the intersection points
 $p_{1}, \ldots, p_{4}$  between $f-2e_{10}, \ldots, f-2e_{13}$ and
 $e_9$, we have a symplectically embedded $2$-sphere, denoted by $S$,
 in $E(1)\sharp 4{\overline{\bf CP}}^2$ which represents a homology class
 $(f-2e_{10})+\cdots +(f-2e_{13})+e_9 =4f +e_9-2(e_{10}+\cdots +e_{13})$
 with square $-9$.
 Now, using a linear plumbing manifold consisting of $5$ disk bundles
 $\{S_{1}, S_{2}, \ldots, S_{5}\}$ lying in a neighborhood of an
 $\widetilde{E_6}$-singular fiber (Figure~\ref{E6}), we obtain a
 configuration $C_{7} \subset E(1)\sharp 4{\overline{\bf CP}}^2$
 by setting $u_1=S_1,\, u_2=S_2,\, \ldots,\, u_5=S_5$ and $u_6=S$
 (Figure~\ref{C7}).  Note that all $2$-spheres $u_{i}$ lying in the
 configuration $C_7$ are symplectically embedded. $\ \ \ \Box$ \\

 {\em Remark 2.} Note that there are other candidates for
 a configuration $C_7$ in $E(1)\sharp 4{\overline{\mathbf CP}}^2$
 by choosing a different linear plumbing manifold lying in a neighborhood
 of an $\widetilde{E_6}$-singular fiber. For example, one may choose
 a linear plumbing manifold consisting of $5$ disk bundles
 $\{S_{1}, S_{2}, S_{3}, S_{6}, S_{7}\}$ in Figure~\ref{E6}
 to get a configuration $C_7$. \\

\begin{theorem}
\label{thm-main1}
 There exists a simply connected symplectic $4$-manifold with
 $b_{2}^+ =1$ and $c_{1}^2 =2$ which is homeomorphic, but not
 diffeomorphic, to ${\mathbf CP}^2 \sharp 7{\overline{{\mathbf CP}}^2}$.
\end{theorem}

 {\em Proof} : {\em Construction} - By Proposition~\ref{pro-C7} above,
 we have a symplectically embedded configuration $C_{7}$ in a
 rational surface $X:=E(1) \sharp 4{\overline{{\mathbf CP}}^2}$.
 Hence we get a new smooth $4$-manifold, denoted by
 $X_{7}=X_{0}\cup_{L(49,-6)}B_{7}$, by rationally blowing down
 along the configuration $C_{7}$ in $X=X_{0}\cup_{L(49,-6)}C_{7}$.
 Furthermore, Theorem~\ref{thm-sy} guarantees the existence of a
 symplectic structure on $X_{7}$.

 {\em Properties of $X_{7}$} - Since a circle representing a generator
 of $\pi_{1}(L(49,-6))$ bounds a disk which is a hemisphere of $S_6$
 lying in $\widetilde{E_6}$-singular fiber, $\pi_{1}(X_{0})=1$.
 Hence the simple connectivity of $X_{7}$ follows from
 Van-Kampen's theorem. Furthermore, it satisfies
 $c_{1}^2(X_{7}) = c_{1}^2(X)+6 =2$, so that it is homeomorphic
 to ${\mathbf CP}^2 \sharp 7{\overline{{\mathbf CP}}^2}$
 due to M. Freedman's classification theorem of simply
 connected closed topological $4$-manifolds.
 It only remains to show that $X_{7}$ is not diffeomorphic to
 ${\mathbf CP}^2 \sharp 7{\overline{{\mathbf CP}}^2}$.
 For this, we first claim that the canonical class $K_7$ on $X_7$
 and the corresponding symplectic $2$-form $\omega_7$ on $X_7$
 satisfy \mbox{$K_7 \cdot[\omega_{7}] > 0$}.
 Then,
 since  $b_{2}^-(X_7) \leq 9$ and \mbox{$(-K_7)\cdot [\omega_{7}] < 0$},
 the Seiberg-Witten invariant $SW_{X_{7}}^{\circ}$ is well defined
 and we have $SW_{X_{7}}^{\circ}(-K_7) =
 SW_{X_{7}}^{\circ,-}(-K_7) = SW_{X_{7}}^{-}(-K_7) = \pm 1$, where
 the last equality follows from Theorem~\ref{T-2}. Note that the
 non-triviality of $SW_{X_{7}}^{\circ}$ means that $X_{7}$ does not
 admit a metric of positive scalar curvature, equivalently, it is not
 diffeomorphic to a rational surface. Hence we are done.

 {\em Proof of $K_7 \cdot[\omega_{7}] > 0$} -
 Note that the canonical class $K$ of
 $X=E(1)\sharp 4{\overline{{\mathbf CP}}^2}$ is represented by
 $K = -3h+(e_{1}+ \cdots + e_{13}) = -[f]+(e_{10}+ \cdots + e_{13})$
 and, by modifying B. Li and T. Li's symplectic cone argument
 in~\cite{ll2}, we may assume that the cohomology class $[\omega]$
 of a symplectic $2$-form $\omega$ on $X$ compatible with
 a canonical class $K$ can be represented by
 $ah - (b_{1}e_{1} + \cdots + b_{13}e_{13})$ for some rational numbers
 $a, b_{1}, \ldots, b_{13}$ satisfying
 $a \geq b_1 \geq b_2 \geq \ldots \geq b_{13} \geq 0$ and
 $3a > b_{1} + \cdots + b_{13}$ (refer to Lemma~\ref{lem-E(1,k)}).
 Now, remembering that $u_{6}=S=4f+e_9 -2(e_{10}+\cdots +e_{13})
 = 12h+e_9 -4(e_{1}+\cdots +e_{9})-2(e_{10}+\cdots +e_{13})$
 in Proposition~\ref{pro-C7} and the canonical class $K$ does not
 intersect with holomorphic $2$-spheres $u_{i}\, (1\leq i \leq 5)$ of
 square $-2$, let us express two cohomology classes $K|_{C_{7}}$ and
 $[\omega|_{C_{7}}]$ using a dual basis $\{\gamma_{i} : 1 \leq i \leq 6\}$
 (i.e. $<\gamma_{i}\,\,,\,u_{j}>=\delta_{ij}$) for $H^{2}(C_{7};{\mathbf Q})$:
\begin{eqnarray*}
 K|_{C_{7}} & = & (K\cdot u_{1})\gamma_{1} + (K\cdot
 u_{2})\gamma_{2}+ \cdots + (K\cdot u_{6})\gamma_{6} \\
 & = & 7\gamma_{6} \ \ \mathrm{and} \\
 \omega|_{C_{7}} & = &
 ([\omega] \cdot u_{1})\gamma_{1} + ([\omega]
         \cdot u_{2})\gamma_{2}+ \cdots + ([\omega] \cdot u_{5})\gamma_{5}
         +([\omega] \cdot u_{6})\gamma_{6} \\
 & = & (b_4-b_7)\gamma_{1} + (b_1-b_4)\gamma_{2} + (a-b_1-b_2-b_3)\gamma_{3}
       + (b_2-b_5)\gamma_{4} \\
 &  &  + (b_5-b_9)\gamma_{5} + \{12a-4(b_1+\cdots +b_9)
       -2(b_{10}+\cdots + b_{13})+b_9 \}\gamma_{6}.
\end{eqnarray*}
 Then, using the intersection form $Q_7$ on $H^{2}(C_{7};{\mathbf Q})$,
 we have
\begin{eqnarray*}
 K|_{C_{7}}\cdot [\omega|_{C_{7}}] & = &
  \frac{-1}{7}\{(b_4-b_7)+ 2(b_1-b_4)+ 3(a-b_1-b_2-b_3)+4(b_2-b_5) \\
 & & \ \ \ \ \ \ + 5(b_5-b_9)+6(12a-4(b_1+\cdots +b_9)
                 -2(b_{10}+\cdots + b_{13})+b_9)\} \\
 & = & \frac{-1}{7}\{75a -25b_1 -23b_2 -27b_3 -25b_4 -23b_5 -24b_6
                     -25b_7 -24b_8 \\
 & & \ \ \ \ \ \ -23b_9 -12(b_{10}+\cdots+ b_{13})\}.
\end{eqnarray*}
 Hence, by Lemma~\ref{lem-can} and Lemma~\ref{lem-E(1,k)}, we have
 \begin{eqnarray*}
   K_7 \cdot [\omega_{7}] & = & K \cdot [\omega]
             - K|_{C_{7}} \cdot [\omega|_{C_{7}}] \\
  & = & \{-3a+(b_{1}+\cdots +b_{13})\} - K|_{C_{7}} \cdot [\omega|_{C_{7}}] \\
  & = & \frac{1}{7}\{54a-18b_1 -16b_2 -20b_3 -18b_4 -16b_5 -17b_6 -18b_7
                     -17b_8 \\
  &   & \ \ \ \  -16b_9 -5(b_{10}+\cdots+ b_{13})\} \\
  & > & \frac{1}{7}\{2b_2 -2b_3 +2b_5 +b_6 +b_8 +2b_9
                     +13(b_{10}+\cdots+ b_{13})\} \\
  & \geq & 0 \, .
 \ \ \ \Box \\
\end{eqnarray*}

 {\em Remark 3.} Since the canonical class $K_7$ induced from a
 symplectic structure $\omega_7$ on the symplectic $4$-manifold $X_7$
 constructed in the proof of Theorem~\ref{thm-main1} above satisfies
 $K_7 \cdot [\omega_{7}] > 0$, one can also conclude directly from
 Lemma~\ref{lem-ll1} that $X_7$ is not diffeomorphic to a rational surface
 ${\mathbf CP}^2 \sharp 7{\overline{{\mathbf CP}}^2}$. \\

 {\em Remark 4.} Similarly, using various different configurations $C_7$
 lying in $E(1)\sharp 4{\overline{\mathbf CP}}^2$ (refer to Remark 2)
 and the same technique as in the proof of Theorem~\ref{thm-main1} above,
 we can construct a family of simply connected symplectic $4$-manifolds
 with $b_{2}^+ =1$ and $c_{1}^2 =2$ which are all homeomorphic,
 but not diffeomorphic, to a rational surface
 ${\mathbf CP}^2 \sharp 7{\overline{{\mathbf CP}}^2}$.
 But we do not know whether all these symplectic $4$-manifolds are
 mutually diffeomorphic to each other. \\

 {\em Remark 5.} There are still some important questions to be solved
 regarding on the symplectic $4$-manifold $X_7$. For example,
 though it is likely to be minimal, it is not easy to determine whether
 $X_7$ is minimal. Furthermore, it is a very intriguing question whether
 $X_7$ admits a complex structure. \\

 {\em Remark 6.} Theorem~\ref{thm-main1} above enables us to confirm that
 a rational surface ${\mathbf CP}^2 \sharp 7{\overline{{\mathbf CP}}^2}$
 admits an exotic smooth structure. \\

\section{Symplectic $4$-Manifolds with $b_{2}^+ =1$ and $c_{1}^2 =1$}
\label{sec-4}
 As we mentioned in the Introduction, the only known simply
 connected symplectic $4$-manifolds with $b_{2}^+ =1$ and
 $c_{1}^2 =1$ are complex surfaces such as a rational surface
 ${\mathbf CP}^2 \sharp 8{\overline{{\mathbf CP}}^2}$ and Barlow
 surfaces. In this section, using the same technique as in the proof
 of Theorem~\ref{thm-main1} above, we construct simply connected
 symplectic $4$-manifolds which are homeomorphic, but not diffeomorphic,
 to a rational surface ${\mathbf CP}^2 \sharp 8{\overline{{\mathbf CP}}^2}$.

\begin{proposition}
\label{pro-C5}
 There exists a configuration $C_5$ in a rational surface
 $E(1)\sharp 3{\overline{\mathbf CP}}^2$ such that all $2$-spheres $u_{i}$
 lying in $C_5$ are symplectically embedded.
\end{proposition}

 {\em Proof} : As in the proof of Proposition~\ref{pro-C7}, we
 blow up $3$ times at the double points of singular fishtail
 fibers so that there exist embedded $2$-spheres,
 $f-2e_{10}, \, \ldots, f-2e_{12}$, in $E(1)\sharp 3{\overline{\bf CP}}^2$
 which intersect a section $e_2$ of $E(1)$ positively at points,
 say $p_1, \, \ldots, p_3$, respectively.
 And then, resolving symplectically the intersection points
 $p_1, \ldots, p_3$  between $f-2e_{10}, \ldots, f-2e_{12}$ and $e_2$,
 we have a symplectically embedded $2$-sphere, denoted by $S$, in
 $E(1)\sharp 3{\overline{\bf CP}}^2$ which represents a homology class
 $(f-2e_{10})+\cdots +(f-2e_{12})+e_2=3f +e_2 -2(e_{10}+\cdots+e_{12})$
 with square $-7$.
 Now, using a linear plumbing manifold consisting of $3$ disk bundles
 $\{S_1, S_2, S_3\}$ lying in a neighborhood of an
 $\widetilde{E_6}$-singular fiber (Figure~\ref{E6}), we obtain a
 configuration $C_5 \subset E(1)\sharp 3{\overline{\bf CP}}^2$
 by setting $u_1=S_1,\, u_2=S_2,\, u_3=S_3$ and $u_4=S$. $\ \ \ \Box$ \\

\begin{theorem}
\label{thm-main2}
 There exists a simply connected symplectic $4$-manifold with
 $b_{2}^+ =1$ and $c_{1}^2 =1$ which is homeomorphic,
 but not diffeomorphic, to a rational surface
 ${\mathbf CP}^2 \sharp 8{\overline{{\mathbf CP}}^2}$.
\end{theorem}

 {\em Proof} : By Proposition~\ref{pro-C5} above, we have a
 symplectically embedded configuration $C_{5}$ in a rational
 surface $X:=E(1) \sharp 3{\overline{{\mathbf CP}}^2}$.
 Hence we get a new smooth $4$-manifold, denoted by
 $X_{5}=X_{0}\cup_{L(25,-4)}B_{5}$, by rationally blowing down
 along the configuration $C_{5}$ in $X=X_{0}\cup_{L(25,-4)}C_{5}$.
 Then the rest of proof is exactly the same as the proof of
 Theorem~\ref{thm-main1} above except a computation of
 $K_5 \cdot [\omega_{5}] > 0$, which is the following:

 In this case, we have $K = -3h+(e_{1}+ \cdots + e_{12})$ and
 $[\omega]= ah - (b_{1}e_{1} + \cdots + b_{12}e_{12})$
 for some rational numbers $a, b_{1}, \ldots, b_{12}$ satisfying
 $a \geq b_1 \geq b_2 \geq \ldots \geq b_{12} \geq 0$ and
 $3a > b_{1} + \cdots + b_{12}$.
 Now, using $u_{4}=S=3f+e_2 -2(e_{10}+\cdots +e_{12})
 = 9h+e_2 -3(e_{1}+\cdots +e_{9})-2(e_{10}+\cdots +e_{12})$,
 we have
\begin{eqnarray*}
 K|_{C_{5}} &=& (K\cdot u_{1})\gamma_{1}+\cdots+(K\cdot u_{4})\gamma_{4}
                =  5\gamma_{4} \ \ \mathrm{and} \\
 \omega|_{C_{5}} & = &
 ([\omega] \cdot u_{1})\gamma_{1}+\cdots + ([\omega] \cdot u_{4})\gamma_{4}\\
 & = & (b_4-b_7)\gamma_{1} + (b_1-b_4)\gamma_{2} +(a-b_1-b_2-b_3)\gamma_{3}\\
 &  &  + \{9a-3(b_1+\cdots +b_9) -2(b_{10}+\cdots + b_{12})+b_2 \}\gamma_{4}.
\end{eqnarray*}
 Hence, by Lemma~\ref{lem-can} and Lemma~\ref{lem-E(1,k)}, we get
\begin{eqnarray*}
 K|_{C_{5}}\cdot [\omega|_{C_{5}}] & = &
     \frac{-1}{5}\{(b_4-b_7)+ 2(b_1-b_4)+ 3(a-b_1-b_2-b_3) \\
 &   & \ \ \ \ \ \ +4(9a-3(b_1+\cdots +b_9)-2(b_{10}+\cdots + b_{12})+b_2)\} \\
 & = & \frac{-1}{5}\{39a-13b_1 -11b_2 -15b_3 -13b_4 -12(b_5+b_6) -13b_7 \\
 &   & \ \ \ \ \ \ -12(b_8+b_9) -8(b_{10}+\cdots+ b_{12})\}, \\
  K_5 \cdot [\omega_{5}] & = & K \cdot [\omega]
             - K|_{C_{5}} \cdot [\omega|_{C_{5}}] \\
 & = & \frac{1}{5}\{24a -8b_1 -6b_2 -10b_3 -8b_4 -7(b_5+b_6) -8b_7 \\
 &   & \ \ \ \      -7(b_8+b_9) -3(b_{10}+\cdots+ b_{12})\} \\
 & > & \frac{1}{5}\{2b_2 -2b_3 +b_5 +b_6 +b_8 +b_9 +5(b_{10}+\cdots+ b_{12})\} \\
 & \geq & 0 \, . \ \ \ \ \ \ \ \Box \\
\end{eqnarray*}

 {\em Remark 7.} Similarly, using various different configurations $C_5$
 lying in $E(1)\sharp 3{\overline{\mathbf CP}}^2$, we can also construct
 a family of simply connected symplectic $4$-manifolds with $b_{2}^+ =1$
 and $c_{1}^2 =1$ which are all homeomorphic, but not diffeomorphic,
 to a rational surface ${\mathbf CP}^2 \sharp 8{\overline{{\mathbf CP}}^2}$.
 But we do not know whether one of these symplectic $4$-manifolds is
 diffeomorphic to a Barlow surface. \\

 Finally, we also investigate exotic smooth structures on a
 rational surface ${\mathbf CP}^2 \sharp 8{\overline{{\mathbf CP}}^2}$.
 As mentioned in the Introduction, Barlow surface, denoted by $B$,
 is homeomorphic, but not diffeomorphic, to a rational surface
 ${\mathbf CP}^2 \sharp 8{\overline{{\mathbf CP}}^2}$.
 Furthermore, whereas ${\mathbf CP}^2 \sharp 8{\overline{{\mathbf CP}}^2}$
 admits an Einstein metric with positive scalar curvature,
 $B$ admits an Einstein metric with negative scalar curvature.
 Now, since a symplectic $4$-manifold $X_7$ constructed
 in Theorem~\ref{thm-main1} has a nontrivial SW-basic class
 $K_{7}$ obtained by a small generic perturbation, so that
 its blow-up manifold $X_7\sharp \overline{{\mathbf CP}}^2$ has
 at least two (up to sign) SW-basic classes,
 $K_{7}\!+\!E$ and $K_{7}\!-\!E$
 ($E$ is an exceptional curve obtained by a blowing up),
 we conclude that the blow-up manifold $X_7\sharp \overline{{\mathbf CP}}^2$
 is also homeomorphic, but not diffeomorphic, to both $B$ and
 ${\mathbf CP}^2 \sharp 8{\overline{{\mathbf CP}}^2}$.
 Furthermore, D. Kotschick pointed out that
 $X_7\sharp \overline{{\mathbf CP}}^2$ does not admit an Einstein
 metric, which can be deduced from his result regarding on
 the existence of Einstein metrics (Theorem~\ref{thm-main3} below
 for details). Hence we get the following useful information about exotic
 smooth structures on a rational surface
 ${\mathbf CP}^2 \sharp 8{\overline{{\mathbf CP}}^2}$.
 Before we state our final result, we first quote a theorem of D. Kotschick:

\begin{theorem}[\cite{k2}]
\label{thm-k2}
  Let $X$ be a smooth $4$-manifold with a monopole class $c$. If
 $X$ admits an Einstein metric, then the maximal number $k$ of
 copies of $\overline{{\mathbf CP}}^2$ that can be split off
 smoothly is bounded by
 \[ k \leq \frac{1}{2}\{2e(X) + 3\sigma(X) - 8d\}\]
 where $d$ is the dimension of moduli space of the solutions of
 Seiberg-Witten equations corresponding to a monopole class $c$.
\end{theorem}

\begin{theorem}
\label{thm-main3}
 A rational surface ${\mathbf CP}^2 \sharp 8{\overline{{\mathbf CP}}^2}$
 admits at least three distinct smooth structures - an Einstein metric
 with positive scalar curvature, an Einstein metric with negative
 scalar curvature and no Einstein metric.
\end{theorem}

 {\em Proof} : Since $X_7\sharp \overline{{\mathbf CP}}^2$,
 $B$ and ${\mathbf CP}^2 \sharp 8{\overline{{\mathbf CP}}^2}$
 are not mutually diffeomorphic to each other, and
 since $B$ and ${\mathbf CP}^2 \sharp 8{\overline{{\mathbf CP}}^2}$
 have desired Einstein metrics, it only suffices to show that
 $X_7\sharp \overline{{\mathbf CP}}^2$ does not admit an Einstein metric.
 For this, suppose that a symplectic $4$-manifold
 $X_7\sharp \overline{{\mathbf CP}}^2$
 has an Einstein metric and apply Theorem~\ref{thm-k2} above to
 $X_7\sharp \overline{{\mathbf CP}}^2$. Then we have
 $1 \leq k \leq \frac{1}{2}\{2e(X) + 3\sigma(X) - 8d\} = \frac{1}{2}$,
 which is a contradiction. $\ \ \ \Box$ \\

\end{document}